\documentclass[10pt,draftclsnofoot,onecolumn]{IEEEtran}  
\ifCLASSINFOpdf
\else
\fi

\usepackage{hyperref}
\usepackage{url}
\usepackage{graphicx}
\usepackage{graphics} 
\usepackage{epsfig} 
\usepackage{mathptmx} 
\usepackage{mathrsfs}
\usepackage{times} 
\usepackage{amsmath} 
\usepackage{amssymb}  
\usepackage{color}
\usepackage{float}
\usepackage{algorithm}
\usepackage{algorithmic}
\usepackage[utf8]{inputenc} 
\usepackage[T1]{fontenc}    
\usepackage{booktabs}       
\usepackage{amsfonts}       
\usepackage{nicefrac}       
\usepackage{microtype}      
\usepackage{color}
\usepackage{amsmath}
\usepackage{amsthm}
\usepackage{enumitem}
\usepackage{wrapfig}
\usepackage{lipsum}
\usepackage{multirow}
\usepackage{booktabs}
\usepackage{array}
\newcolumntype{P}[1]{>{\centering\arraybackslash}p{#1}}

\newcommand{\R}{\mathbb{R}}

\newtheorem{theorem}{Theorem}

\newtheorem*{lemma*}{Lemma}
\newtheorem{definition}{Definition}
\newcommand{\bd}{\mathbf}
\newcommand\ChangeRT[1]{\noalign{\hrule height #1}}

\newtheorem{Definition}[theorem]{\textbf{Definition}}
\newtheorem{Proposition}[theorem]{\textbf{Proposition}}

\DeclareMathAlphabet{\mathpzc}{OT1}{pzc}{m}{it}

\begin{document}
\title{Input Convex Neural Networks for Optimal Voltage Regulation}

\author{Yize Chen,
        Yuanyuan Shi, and
        Baosen Zhang
\thanks{The authors are with the Department of Electrical and Computer Engineering at the University of Washington, emails:\{yizechen,yyshi,zhangbao\}@uw.edu.}

}


\maketitle


\begin{abstract}
The increasing penetration of renewables in distribution networks calls for faster and more advanced voltage regulation strategies. A promising approach is to formulate the problem as an optimization problem, where the optimal reactive power injection from inverters are calculated to maintain the voltages while satisfying power network constraints. However, existing optimization algorithms require the exact topology and line parameters of underlying distribution system, which are not known for most cases and are difficult to infer. In this paper, we propose to use specifically designed neural network to tackle the learning and optimization problem together. In the training stage, the proposed input convex neural network learns the mapping between the power injections and the voltages. In the voltage regulation stage, such trained network can find the optimal reactive power injections by design. We also provide a practical distributed algorithm by using the trained neural network. Theoretical bounds on the representation performance and learning efficiency of proposed model are also discussed. Numerical simulations on multiple test systems are conducted to illustrate the performance of the algorithm.

\end{abstract}


\section{Introduction}
\label{Sec:intro}
The power distribution system is current undergoing a series of rapid transformations. With the growing adoption of distributed energy resources~(e.g., rooftop PV and electric vehicles), distribution systems are experiencing greater variations in the active power injections~\cite{pedrasa2010coordinated}. Because of the high $r/x$ ratios of transmission lines in the distribution system, the voltage magnitude at the buses are sensitive to active power injections and larger variations in the output of distributed energy resources (DERs) often lead to unacceptable voltage swings~\cite{carvalho2008distributed,tonkoski2010coordinated,keane2010enhanced}. Existing voltage regulation devices such as tap-changing transformers and switched capacitors are effective in dealing with \emph{slow} variations in voltage~(on the timescale of hours), but not variations induced by the DERs~(in the timescale of minutes)~\cite{joos2000potential}. Because of their mechanical nature, using these devices for fast voltage regulation may significantly degrade their service lifetimes~\cite{logue2001utility,kersting2006distribution}.


\begin{figure}[ht]
	\centering
	\includegraphics[width=0.7 \columnwidth]{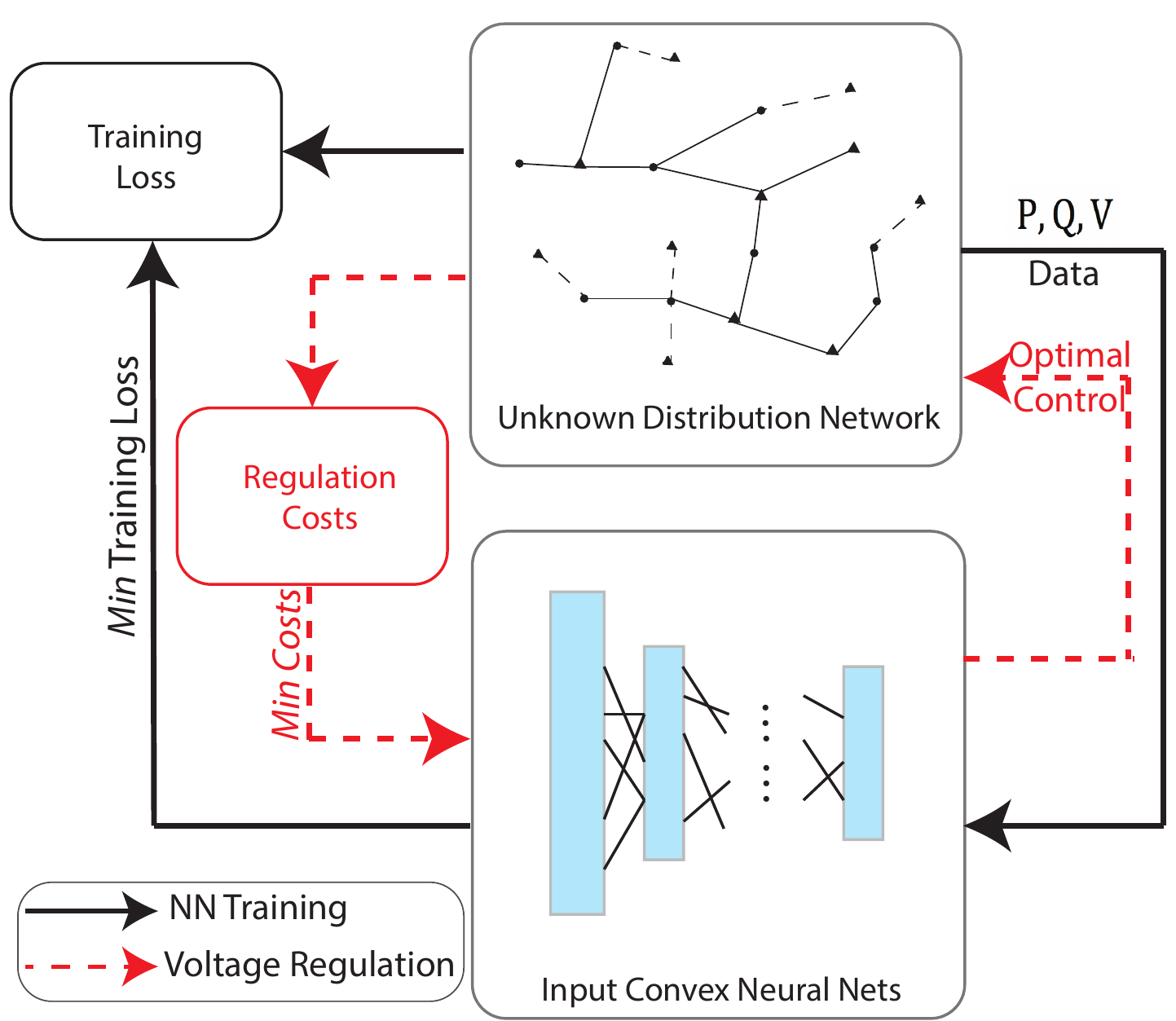}
	\caption{\small Proposed data-driven method for distribution networks with unknown topology. An input convex neural network is fitted to learn the mapping from power injections to a cost function based on voltage, and then an optimization problem is solved to find the optimal reactive power injections.}
	\label{fig:schematic}
\end{figure}

Instead of using dedicated mechanical devices, the power electronic interfaces of the DERs themselves can be used as a mean to control active and reactive power injections for voltage regulations. Many current inverters are endowed with wireless or wireline communication capabilities and can set their operating points electronically. A growing body of literature have focuses on this approach, most of which focuses on casting the voltage regulation problem as an optimization problem, then proposing various algorithms in both centralized and decentralized settings. For a more complete review, please see~\cite{madani2015promises,molzahn2017survey,zhu2015fast} and the references within.

Most relevant to this paper is the approach of formulating the voltage regulation problem as an optimal power flow (OPF) problem. A line of work has showed that OPF problems can be made convex when the network has a tree topology~(see, e.g.~\cite{low2014convex,Zhang13_geometry} and the references within), and since most distribution systems operate as a tree (i.e., radially)~\cite{kersting2006distribution}, the voltage regulation problem is in fact convex under many settings~\cite{farivar2012optimal,zhang2014optimal}. More specifically, the OPF problem for a tree network can be shown to have exact SOCP or SDP convex relaxations~\cite{boyd2004convex}. When the network is large or computational resources are limited, linearized or other approximate models have been proposed to solve centralized or distribution versions of voltage regulation problems~\cite{qu2019optimal,magnusson2019voltage, tang2018fast}.

A common basis of the optimization based approach for voltage regulation is the knowledge of the topology and the parameter of the physical distribution system. However, for many distributions systems, neither the topology nor the line parameters are known by the system operators~\cite{abur2004power,weng2013graphical,li2019robust}. Unlike transmission systems, distribution systems are typically not monitored, while the (secondary) network topology is rarely digitized. Besides, learning the topology and the line parameters (the $\bd Y$-bus matrix) is a hard problem. Existing learning algorithms either require full PMU measurements~\cite{liao2018unbalanced} or a large number of time-stamped smart meter measurements~\cite{luan2015smart}. The former is rarely available for a distribution network, and the latter may take too long.

 In this work, we combine the availability of smart meter data with the existing results on the convexity of the OPF problem to enable optimal voltage regulation for \emph{distribution networks with unknown topology and parameters}. Namely, we use a structured neural network surrogate model of the distribution system that is guaranteed to be convex from the power injections to the bus voltage magnitudes. We call this network an \emph{input convex neural network~(ICNN)} \cite{amos2017input,chen2018optimal}, which is constructed using common neural network activation functions with constrained weights and can be trained using smart meter data and standard back propagation algorithms. The precise construction is given later in the paper, but once it is trained, it is then used in a model predictive manner to solve the voltage regulation problem. As shown in Fig. \ref{fig:schematic}, in the training stage, ICNN is using measurement data (we assume the active, reactive and voltage magnitude information are available, while the phase is not measured) to learn the unknown, nonlinear mappings from active and reactive power injections to bus voltages; in the voltage regulation stage, ICNN serves as the model to be optimized over.

 Our approach is along the classical two-step approach of first performing system identification then optimizing over the identified model. The key insight is to strike a tradeoff between these two steps. Since without phase information, directly learning the $\bd Y$-bus matrix is difficult~\cite{liao2018unbalanced}, the most accurate model is to fit a deep neural network (DNN) using data on power injections and voltage magnitudes. Although accurate, normal-trained deep learning model is intractable to use in the following optimization stage since DNNs are typically highly nonconvex from input to output and can have a large number of local optima~\cite{chen2018optimal}. Therefore, a linear model is often used in practice to ensure the optimization problem can be solved efficiently~\cite{li2014real}. This approach, however, suffers from the fact that the mapping from power injection to voltage magnitude is not linear, especially for heavily loaded systems~\cite{Zhang2013}. In this paper, we show that the ICNN achieves the optimal tradeoff between computational tractability and model accuracy:  it captures the nonlinear mapping accurately (since we know the physical model is convex) and ensuring the optimization problem is convex.

 There are other lines of research on voltage regulation that do not follow the above two-step approach. For example, a completely local controller can be designed (e.g., voltage droop~\cite{guo2014distributed}) that injects reactive power based on the local voltage measurements~\cite{Zhang2013}. It is known that linear control laws are insufficient, in the sense that there always exist systems that cannot be stabilized by linear laws~\cite{Zhang2013}, or experience  efficiency issues when the network is large~\cite{zhu2015fast}. More recently, machine learning techniques have been used to search for nonlinear control laws~\cite{jalali2019designing}. In another approach, reinforcement learning is proposed to learn a (centralized or decentralized) policy for voltage regulation~\cite{yang2019real}. Although this approach either requires a detailed model of the system for simulation or an actual system that can be repeatedly experimented with. The former is not possible if the system is unknown, while the latter may be impractical in real systems design and implementations.

 In summary, we make the following contributions to the voltage regulation problem:
 \begin{enumerate}
 	\item We consider distribution systems with unknown topology and parameters;
 	\item We design centralized and decentralized data-driven controller with optimality guarantee, providing an alternative to the difficult problem of explicitly finding the $\bd Y$-bus matrix of the network;
 	\item We prove the representation power and efficiency of ICNNs, and show how it leverages existing results about the convexity of the problem.
\end{enumerate}

 The remainder of the paper is organized as follows.  Section \ref{Sec:Model} discusses the voltage regulation problem in distribution networks and formulates it as an optimization problem. Section \ref{Sec:ICNN} states the main learning and control framework, and provides theoretical guarantees on the representation power, learning efficiency as well as control optimality. Section~\ref{sec:voltage} shows how the learned neural network can be used for voltage regulation and Section \ref{Sec:Distributed} extends the proposed algorithm to the decentralize setting. The performances are illustrated in Section \ref{Sec:Experiment} via numerical experiments on IEEE standard test cases. Concluding remarks are presented in Section \ref{Sec:Conclusion}.

\section{Model and Problem Formulation}
\label{Sec:Model}
\subsection{Power Flow Equations}
In this paper we consider distribution systems with a tree topology, with a set $\mathpzc{N}=\{1,...,N\}$ of buses and a set $\mathpzc{E}\in \mathpzc{N} \times \mathpzc{N}$ of lines. For each bus $i$, denote $V_i$ as the voltage magnitude and $\theta_{i}$ as the voltage phase angle; let $p_i$ and $q_i$ denote the active and reactive power injections; let $s_i=p_i+jq_i$ be the complex power injection at bus $i$. The corresponding active and reactive power injection vectors are denoted as $\mathbf{p}=\left[\begin{array}{llll}{p_{1}} & {p_{2}} & {\cdots} & p_{N}\end{array}\right]^{T}, \mathbf{q}=\left[\begin{array}{llll}{q_{1}} & {q_{2}} & {\cdots} & {q_{N}}\end{array}\right]^{T}$.
For each line $(i,k)\in \mathpzc{E}$, denote line admittance $y_{ik}=g_{ik}-j b_{ik}$ with $b_{ik}>0,\; g_{ik}>0$.

The power flow equations can be stated in many forms, and we use the DistFlow formulation~\cite{baran1989optimal} for the ease of comparison with SOCP relaxations and linear approximations. In this formulation, let $s_{ik}=p_{ik}+jq_{ik}$ denote the complex power flow from bus $i$ to bus $k$ and $z_{ik}=1/y_{ik}=r_{ik}+jx_{ik}$ denote the line impedance. The power flow equations are
\begin{subequations}
	\label{equ:branch_flow}
	\begin{align}
-p_{k} &=p_{i k}-r_{i k} l_{i k}-\sum_{l :(k,l) \in \mathpzc{E}} p_{kl} \label{subequ:branch_flowP}\\
-q_{k} &=p_{i k}-x_{i k} l_{i k}-\sum_{l :(k,l) \in \mathpzc{E}} q_{kl} \label{subequ:branch_flowQ}\\
V_{k}^2 &=V_{i}^2-2\left(r_{i k} p_{i k}+x_{i k} q_{i k}\right)+\left(r_{i k}^{2}+x_{i k}^{2}\right)l_{i k} \label{subequ:branch_flowV}\\
l_{i k} &=\frac{p_{i k}^{2}+q_{i k}^{2}}{V_{i}^2}.\label{subequ:branch_flowl}
		\end{align}
	\end{subequations}

The optimal voltage regulation problem is then to optimize system performance subject to the power flow equations \eqref{equ:branch_flow}. There are a variety of objective functions and constraints we can consider~\cite{zhu2015fast}, and in this paper we adopt the following optimal voltage regulation formulation:
\begin{subequations}
	\label{eqn:voltage_regulation}
	\begin{align}
     \min_{\mathbf{q}} \quad &\sum_{i=1}^N \alpha_i |V_i-V_{i,0}| \label{subequ:loss}\\
	\text { s.t. } \quad &\underline {\mathbf{q} }\leq \mathbf{q} \leq \overline{\mathbf{q}} \label{eqn:constraint_reactive}\\
	& \text{ Power Flow Equations } \eqref{equ:branch_flow}  \label{eqn:constraint_flow}
	\end{align}
\end{subequations}
where $\alpha_i$ is a weight parameter. The goal of the problem is to maintain voltage magnitude $V_i$ within a small distance from the nominal value $V_{i,0}$ for all buses~(e.g., plus/minus $5\%$)~\cite{li2018distribution}. The constraints on reactive powers comes from the rating of the power electronics on the DERs. Note that if a bus does not have reactive power capability, we can simply set the upper and lower bound to be equal to the nominal reactive load value.
Here we assume that the active power injections are fixed exogenously (e.g., by solar irradiation), although we can accommodate active power as a optimization variable using the same methodology presented in this paper. In addition, other cost functions such as losses or costs or reactive power can be added to the objective as well.

The optimization problem in \eqref{eqn:voltage_regulation} is not convex because of the quadratic equality~\eqref{subequ:branch_flowl} in the power flow equations. A simple relaxation is to make it into an inequality as:
 \begin{equation}
 \label{equ:branch_flow_relaxed}
l_{i k} \geq \frac{p_{i k}^{2}+q_{i k}^{2}}{V_{i}^2}.
 \end{equation}
 Using \eqref{equ:branch_flow_relaxed} instead of \eqref{subequ:branch_flowl} in \eqref{eqn:voltage_regulation} makes it a convex optimization problem, in particular a second order cone problem~(SOCP)~\cite{farivar2012optimal,farivar2011inverter}.

 As outlined in the introduction, an important result is that the above convex relaxation is tight under many settings. In fact, there are different convex relaxations one can use, including semidefinite programming (SDP) and the result extends to unbalanced systems as well~\cite{peng2015distributed,watson2017optimized}. Therefore, in principle, there exist efficient algorithms that are guaranteed to find the global optimal solution. However, as we discuss in the next section, these algorithms have rarely been implemented in practice~\cite{liao2018unbalanced}.

\subsection{Practical Design Requirements}
A major challenge in solving the voltage regulation problem (or other optimization problem in the distribution system) is that the topology and linear parameters of the system are often not known. Due to the lack of monitoring, poor digital record keeping, age of the networks, frequent repairs and switching, utilities may not have a digitized model of their systems. Even as the amount of sensors have grown exponentially in the system, most of these are smart meters (AMI sensors) and do not have synchronous fine-grained voltage phase measurements. Unfortunately, without phase information, it seems that directly learning the $\bd Y$-bus matrix is difficult~\cite{deka2016estimating,li2019robust}.

 Given these practical challenges of voltage regulation, we want to design a strategy that satisfies following requirements:
 \begin{enumerate}
 	\item The controller must learn an accurate representation of the power injections to nodal voltage magnitudes using data without voltage phase measurements;
 	\item Such representation is easily integrated into the optimization framework.
\end{enumerate}
 Intuitively, we are trying to design and find functions $\left| V_i-V_{i,0}\right|=f_{i}(\mathbf{p,q}), i=1,...,N$, which could accurately represent the relationship from active and reactive power injections to nodal voltage magnitude deviation. By leveraging historical smart meter data to fit $f_i$, we want the fitted model to faithfully represent the underlying grid. Importantly, if $f_i$ is a convex function from $\mathbf{p,q}$ to $\left| V_i-V_{i,0}\right|$, then the following problem
 \begin{subequations}
 	\label{eqn:main2}
 	\begin{align}
 	\min_{\mathbf{q}} \quad & \sum_{i=1}^N \alpha_i |V_i-V_{i,0}| \label{subeqn:obj}\\
 	\text { s.t. } \quad &\underline {\mathbf{q} }\leq \mathbf{q} \leq \overline{\mathbf{q}} \label{eqn:constraints2}\\
 	& \left| V_i-V_{i,0}\right| = f_{i}(\mathbf{p,q})\label{subeqn:f}
 	\end{align}
 \end{subequations}
 is still a tractable convex optimization problem. Since we know convex relaxations of the power flow equations are tight, restricting $f_i$ to be convex does not result in a loss of generality.  In the next section, we will describe how we design the learning model $f_i$ based on neural networks, which not only efficiently learns the mapping from active and reactive power injections to the voltage magnitude deviations more powerfully than a linear model, but is also guaranteed to be a convex function.

\section{Input Convex Neural Network Design}
\label{Sec:ICNN}
In this section, starting from the standard neural networks architecture, we illustrate how to construct a neural network whose outputs are convex with respect to inputs. We then show how to apply such input convex neural networks (ICNN) in the task of voltage regulation in distribution networks, and describe a practical algorithm to find optimal reactive power injections under reactive power capacity constraints.

\subsection{Neural Networks for Function Fitting}
For a standard setup of neural networks~(NN) model,  the multi-layer network is composed of an input layer $x$, $m$ hidden layers $z_l,\;l=1,...,m$ with parameters $\theta:=\{W_l,b_l\} \;i=1,...,m$, and an output $y$. For notation simplicity, we use $h_\theta(x)$ to denote the neural networks with input $x$ and parameters $\theta$. For the computation at layer $l$, an nonlinear activation function $g_i( \cdot)$ is used. For instance, rectified linear unit (ReLU) is a popular choice with $g(x)=\max(0,x)$. Given input $x$, a neural network is implementing the following computation
\begin{subequations}
	\begin{align*}
	z_1=&g_1(W_1 x+b_1);\\
	z_{l}=&g_{l+1}(W_{l}z_{l-1}+b_{l}),\quad l=2,...,m
	\end{align*}
\end{subequations}
and the neural network output is the value of the last layer $z_m$. In the task of supervised learning, back-propagation algorithms based on gradient descent are used to train a group of $\{W_l,b_l\}$ that minimize the training loss defined as $L(y, h_{\theta}(x))$~\cite{bottou2010large}. The choice of training loss is task specific and we use the squared loss.
We design specific neural networks such that the output of $h_\theta(x)$ is convex with respect to $x$, which could be then representing $f_i$ in \eqref{subeqn:f} and be used for solving \eqref{eqn:main2} once the model is trained.

\subsection{ICNN Architecture Design}
We adapt the neural networks design from our previous work~\cite{chen2018optimal} and the original input convex neural networks~(ICNN) proposed in \cite{amos2017input} to the setting of multiple inputs (e.g., $\{\mathbf{p,q}\}$) and multiple outputs (e.g., $|V_i-V_{i,0}|, \forall i$). Then by convexity, we mean each of ICNN's output is convex with respect to all the dimension of inputs.

\begin{figure}[ht]
	\centering
	\includegraphics[width=0.65 \columnwidth]{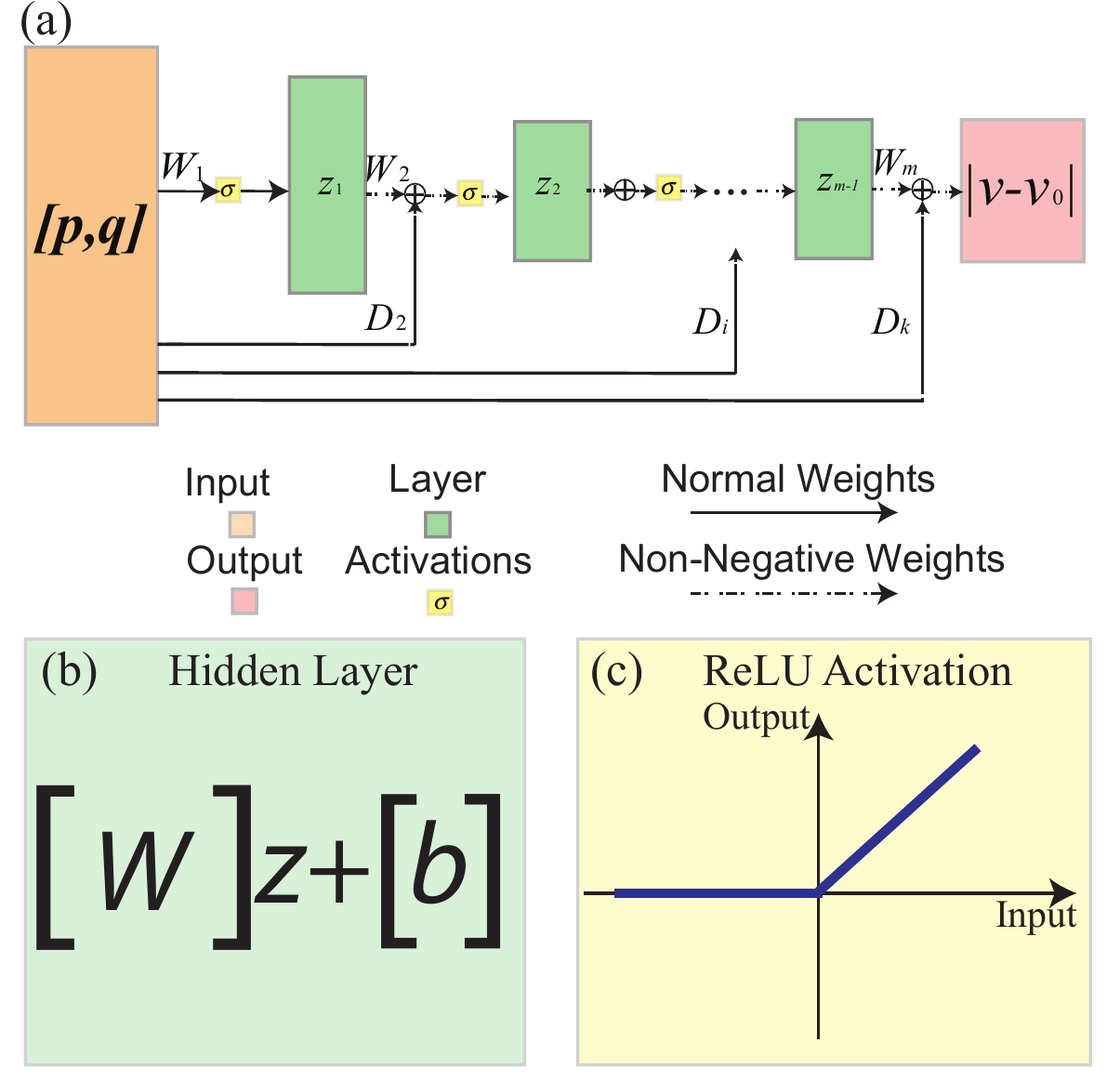}
	\caption{\small (a). The input convex neural networks (ICNN) architecture design; (b) the layer computation for ICNN, where we constrain $W_{2:m}$ to be non-negative; (c). a ReLU nonlinear activation function.}
	\label{fig:NN}
\end{figure}

The following proposition summarizes the major adaptations we make to the standard neural networks that guarantee modeling convexity:
\begin{Proposition}
	\label{prop:icnn_arch}
	The network shown in Fig. \ref{fig:NN}(a) is a convex function from inputs to outputs provided that all $W_{2:m}$ are non-negative, and all $g_l$ are convex and non-decreasing functions.
\end{Proposition}
The restriction on the $g_l$ function to be convex and non-decreasing function is actually not a strong restriction. Popular activation functions like ReLU function shown in Fig.~\ref{fig:NN}(c) already satisifes such restriction. Thus with the nonnegative constraints on $W_{2:m}$ and with the choice of activation functions, we are already constructing a neural network whose output is convex with respect to inputs.

To compensate the neural network representability loss due to the constraints of $W_l \geq 0,\; l=2,...,m$, we add direct passthrough layers $D_l, \; l=2,...m$ from input to subsequent layers, and there is no constraint on the weights of these links. Such direct links have also been widely used in the design of deeper neural networks, which have achieved better performance in various learning tasks~\cite{he2016deep}. Combined with layer bias $b_l$, layer $i$ passes its value through ReLU activation function and goes to next layer $i+1$.

Mathematically, for each layer $l=1,...,m$, the layer-wise computations are modified as follows
\begin{subequations}
	\begin{align}
	z_1=&g_1(W_1 x+b_1);	\\
	z_{l}=&g_{l}\left(W_{l} z_{l-1}+D_{l} x+b_{l}\right), \quad l=2,...,m
	\end{align}
\end{subequations}

We note that our neural networks design naturally extends to the scenario when input $x$ and output $y$ are high-dimensional vectors. Meanwhile, each dimension of output is convex with respect to all the inputs. Such property also makes it possible to fit a single ICNN $h_\theta$ to model multiple convex functions, which will be shown to be useful for voltage regulation task.

The convexity of the proposed neural network directly follows from the composition rule of convex functions~\cite{boyd2004convex}, which states that the composition of an inner convex function and an outer convex, non-decreasing function is convex. The structure of the input convex neural network (ICNN) structure in Proposition~\ref{prop:icnn_arch} is motivated by the structure in~\cite{amos2017input} but modified to be more suitable to control of dynamical systems. In ~\cite{amos2017input} it only requires $\bd W_{2:k}$ to be non-negative while having no restrictions on weights $\bd W_1$ and $\bd D_{2:k}$. Our construction achieves the exact representation by expanding the inputs to include both $\bd u$ ($\in \R^{d}$) and $-\bd u$.


Although it allows for any increasing convex activation functions, in this paper we work with the popular ReLU activation function and its variants. Two notable additions in ICNN compared with conventional feedforward neural networks are: 1) Addition of the direct \emph{``passthrough'' layers} connecting inputs to hidden layers and conventional \emph{feedforward layers} connecting hidden layers for better representation power; \newline 2) the expanded inputs that include both $\bd u$ and $-\bd u$. Note that such construction guarantees that the network is convex and non-decreasing with respect to the expanded inputs $\hat{\bd u}=\begin{bmatrix}\bd u \\-\bd u\end{bmatrix}$, while the output can achieve either decreasing or non-decreasing functions over $\bd u$.
\begin{figure}[ht]
	\centering
	\includegraphics[width=0.65 \columnwidth]{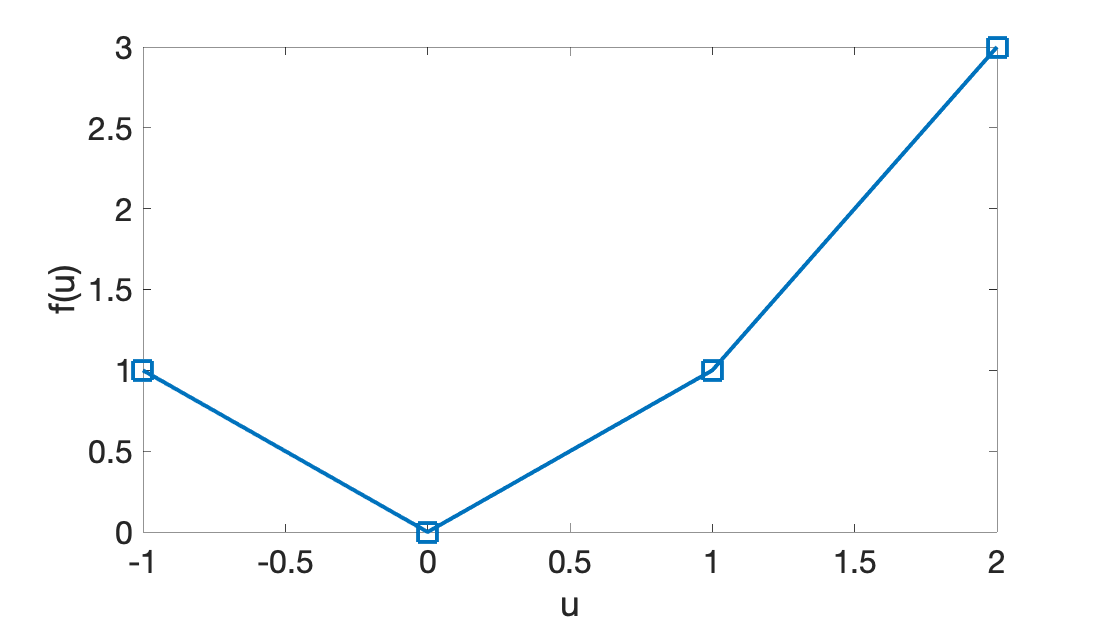}
	\caption{Example of a piecewise convex function that is partly decreasing and partly increasing.}
	\label{fig:example}
\end{figure}

\noindent \textbf{Example} Here we show how to use a ICNN to fit the convex function in Fig.~\ref{fig:example}. The function $f(u)$ in the example has domain $[-1,2]$ and
\begin{equation*}
	f(u)= \begin{cases}
	-u &\mbox { if } -1\leq u \leq 0 \\
	u &\mbox{ if } 0 < u \leq 1 \\
	2u-1 &\mbox{ if } 1 < u \leq 2
\end{cases}.
\end{equation*}
It is easy to check that the following function fits $f(u)$ exactly:
\begin{equation*}
	\begin{bmatrix}1 \\ 1 \\ 1 \end{bmatrix}^T g\left(\begin{bmatrix} 1 & 0 \\ 1 & 0 \\ 0 & 1\end{bmatrix} \begin{bmatrix} u \\ -u \end{bmatrix}+\begin{bmatrix} 0 \\ -1 \\ 0 \end{bmatrix}\right),
\end{equation*}
where $g$ is the componentwise ReLU activation function. Such matrix weights and biases can be also easily learned by ICNN hidder layer.

\subsection{ICNN Representation Power and Efficiency}
Besides the computational traceability of the input convex neural network, a natural question is related to the function approximation capability of the proposed ICNN. For instance, once ICNN is used in \eqref{eqn:main2} for voltage regulation, it should be able to fit the underlying convex functions from $\mathbf{p,q}$  to $\left| V_i-V_{i,0}\right|$ as accurate as possible.  This subsection provides theoretical guarantees that ICNN is able to fit any convex functions, which opens the door to integrate ICNN in convex optimization.

\begin{definition}
	Given a function $f: \mathbb{R}^d \rightarrow \mathbb{R}$, we say that the function $\hat{f}$ approximate $f$ within $\epsilon$ if $|f(\bd x)-\hat{f} (\bd x)| \leq \epsilon$ for all $\bd x$ in the domain of $f$.
\end{definition}
\begin{theorem}\label{thm:power}[Representation power of ICNN]
	For any Lipschitz convex function over a compact domain, there exists a neural network with nonnegative weights and ReLU activation functions that approximates it within $\epsilon$.
\end{theorem}
The proof of this theorem can be found in~\cite{chen2018optimal}.

It turns out to prove Theorem~\ref{thm:power}, we first approximate a convex function by a maximum of affine functions. Then we are able to construct an input convex neural network according to this maximum.  A natural question arises, as why there is a need to learn a neural network rather than directly fitting the affine functions in the maximum? This approach was taken in~\cite{magnani2009convex}, where a convex piecewise-linear function (max of affine functions) are directly learned from data through a regression problem.

A key reason that we propose to use ICNN to fit a function rather than directly finding a maximum of affine functions is that the former is a much more efficient parameterization than the latter. As stated in Theorem~\ref{thm:complexity}, a maximum of $K$ affine functions can be represented by an ICNN with $K$ layers, where each layer only requires a single ReLU activation function. However, given a single layer ICNN with $K$ ReLU activation functions, it may take a maximum of $2^K$ affine functions to represent it exactly. Therefore in practice, it would be much easier to train a good ICNN than finding a good set of affine functions.
%
\begin{theorem}\label{thm:complexity}[Efficiency of Representation]
	\begin{enumerate}[noitemsep,nolistsep]
		\item Let $f_{ICNN}: \mathbb{R}^d \rightarrow \mathbb{R}$ be an input convex neural network with $K$ ReLU activation functions. Then $\Omega(2^K)$ functions are required to represent $f_{ICNN}$ using a max of affine functions.
		\item 	Let $f_{CPL}: \mathbb{R}^d \rightarrow \mathbb{R}$ be a max of $K$ affine functions. Then $O(K)$ activation functions are sufficient to represent $f_{CPL}$ exactly with an ICNN.
	\end{enumerate}
\end{theorem}
The proof of this theorem is again given in~\cite{chen2018optimal}.

\section{ICNN for Voltage Regulation}
\label{sec:voltage}
\begin{algorithm}
	\caption{ICNN for Voltage Regulation}
	\label{algorithm1}
	\begin{algorithmic}
		\REQUIRE Learning rate $\eta$, Step size $\gamma$, Batch size $T$, Training iterations $n_{training}$, Optimization stopping critertion $\epsilon$
		\ENSURE Training dataset $\{\mathbf{p,q,V}\}$
		\ENSURE Initial model  $h_\theta$  \\
		\STATE \emph{\# ICNN Training}
		\FOR {$iter=0,...,n_{training}$}
		\STATE \emph{\# Update parameters for ICNN}
		\STATE Sample batch from historical data:\\ $\{\mathbf{p}^{k},\mathbf{q}^{k},\mathbf{V}^{k}\}_{k=1}^T \sim \mathbb{P}_{x}$
		\STATE Update $h_{\theta}$ using stochastic gradient descent:
		\STATE $\mathbf{V}_{target}^k=\{|V_i^k-V_{i,0}^k|\}, i=1,...,N$
		\STATE $h_{\theta}= h_{\theta}-\eta \Pi_{W_{2:m}\geq0}\left(h_{\theta}-\nabla_{h_\theta}(L(\mathbf{V}_{target}), h_\theta (\mathbf{p},\mathbf{q}))\right)$
		\ENDFOR
		\STATE Fix ICNN parameters $h_\theta$
		\STATE \emph{\# Voltage Regulation via ICNN}
		\STATE Get measurements $\{\mathbf{p,q}\}$, $t=0, \; \mathbf{q}^{(t)}=\mathbf{q}$
		\WHILE{$\nabla_{\mathbf{q}^{(t)}}(\sum_{i=1}^N h_\theta(\mathbf{p,q}^{(t)}))>\epsilon$}
		\STATE $\mathbf{q}^{(t+1)}=\mathbf{q}^{(t)}-\gamma \Pi_{\mathbf{q}}(\mathbf{q}^{(t)}-\nabla_{\mathbf{q}^{(t)}}(\sum_{i=1}^N h_\theta(\mathbf{p,q}^{(t)})))$
		\STATE $t\leftarrow t+1$
		\ENDWHILE
		\STATE Optimal reactive power injection: $\mathbf{q}^*=\mathbf{q}^{(t)}$
	\end{algorithmic}
\end{algorithm}

In this section, we illustrate how to use ICNN for voltage regulation when the underlying topology and the line parameters are unknown. we propose to first learn a convex mapping from $\{\mathbf{p,q}\}$ to voltage magnitude deviations using an ICNN. Once fitted using collected observations, we are able to use the same ICNN, and integrate it to \eqref{eqn:main2}  to find optimal reactive power injections.

\subsection{ICNN Training}
In order to train ICNN and learn its  parameters $h_\theta$, we need to minimize the supervised training loss defined on training data pairs. For the $k$th training instance, it is defined as the mean square error between the ground truth voltage magnitude deviation vector $\mathbf{V}_{target}:=\{|V_i^k-V_{i,0}^k|\}, i=1,...,N$ and the ICNN output:
\begin{equation}
L(\mathbf{V}_{target}, h_\theta(\mathbf{p,q}))=\frac{1}{N}||\mathbf{V}_{target}-h_\theta(\mathbf{p,q})||_2^2,
\end{equation}
and the update of $h_\theta$ is based on gradient descent algorithm. In addition, to take the constraints of $W_{2:m}\geq 0$ into account, we need to make sure the gradient descent update always falls into the feasible regions (e.g., nonnegative weights). Hence we use a projected gradient algorithm to guarantee the constraint holds~\cite{boyd2004convex}.

\begin{Definition}
	The projection of a point $y$, onto a set $X$ is defined as
	\begin{equation}
	\label{equ:Projection}
	\Pi_{X}(y)=\underset{x \in X}{\operatorname{argmin}} \frac{1}{2}\|x-y\|_{2}^{2}
	\end{equation}
\end{Definition}

Given a starting point $x^{(0)}\in X$ and step-size $\gamma>0$, projected gradient descent~(PGD) extends the standard gradient descent settings with the projection step onto the feasible sets of feasible reactive power. At iteration $t$, the algorithm takes the following PGD step:
\begin{equation}
\label{equ:PGD}
x^{(t+1)}=x^{(t)}-\gamma \Pi_{X}\left(x^{(t)}-\nabla f\left(x^{(t)}\right)\right), \forall t \geq 1
\end{equation}
which is implemented iteratively until a certain stoppping criterion~(e.g., fixed number of iterations or gradient value is smaller than predefined $\epsilon$) is satisfied.  The ICNN weights are then updated as follows\vspace{-10pt}

\begin{equation}
\label{eqn:GD}
h_{\theta}= h_{\theta}-\gamma \Pi_{W_{2:m}\geq0}\left(h_{\theta}-\nabla_{h_\theta}(L(\mathbf{V}_{target}), h_\theta (\mathbf{p},\mathbf{q}))\right).
\end{equation}

In practical implementations where there are large groups of measurements $\{\mathbf{p}^{k},\mathbf{q}^{k},\mathbf{V}^{k}\}$ with $k$ standing for the index for measurement index, it is possible to use small batch of training data to do PGD steps \eqref{eqn:GD}. Such practical algorithms, e.g., stochastic gradient descent, can accelerate training convergence \cite{bottou2010large}. The training procedure is summarized in Algorithm \ref{algorithm1} by using collected training data and stochastic gradient descent training algorithm.

\subsection{Integrate ICNN as Optimal Controller}
Once the ICNN training process is finished, we fix model parameters $h_\theta$, and use it as a proxy model for the unknown distribution networks model $f_i,\; i=1,...,N$ in \eqref{subeqn:f}. Since $h_\theta$ represents the convex mappings from $\mathbf{q}$ to $|V_i-V_{i,0}|,\;\forall i$, we are now ready to solve \eqref{eqn:main2} computationally. In the similar spirit of ICNN training, where we optimize over neural network weights using gradient descent to minimize training loss, in the voltage regulation setting, we optimize over ICNN inputs $\mathbf{q}$ to minimize the optimization objective $ \sum_{i=1}^N \alpha_i |V_i-V_{i,0}|$ using gradient descent. Such optimization procedure over ICNN is also highly tractable using standard off-the-shelf machine learning packages such as Tensorflow~\cite{abadi2016tensorflow}. Again, to take the constraints of reactive power injection range into account, we need to make sure that the gradient descent update always falls into the feasible reactive power injection regions.  By adapting PGD to the trained ICNN, starting from uncontrolled reactive power $\mathbf{q}^{(0)}=\mathbf{q}$, we take iterative PGD steps on the voltage regulation objective \eqref{subeqn:obj} until gradient convergence. PGD steps also guarantee the convergence to optimal solution under the convex settings. The overall algorithm for ICNN training and finding optimal reactive power injections are described in Algorithm \ref{algorithm1}.

Fundamentally, ICNN allows us to use neural networks in decision making processes by guaranteeing the solution is unique and globally optimal. Since many complex input and output relationships can be learned through deep neural networks, it is natural to consider using the data-driven, learned model in an optimization problem in the form of
\begin{subequations} \label{eqn:opt_single}
	\begin{align}
	\min_{\bd u} &\; f(\bd u; \bd W) \\
	\mbox{s.t. } &\; \bd u \in \mathpzc{U},
	\end{align}
\end{subequations}
where $\mathpzc{U}$ is a convex feasible space. Then if $f$ is an ICNN, optimizing over $\bd u$ is a convex problem, which can be solved efficiently to global optimality. Note that we will always duplicate the variables by introducing $\bd v =-\bd u$, but again this does not change the convexity of the problem. Of course, since the weights of the network are restricted to be nonnegative, the performance of the network (e.g., classification) may be worse. A common thread we observe in this paper is that trading off neural network fitting performance with tractability can be preferable.

\section{Distributed Algorithm}
\label{Sec:Distributed}
In Section \ref{sec:voltage}, we show that by utilizing a trained ICNN, optimal reactive power injections can be calculated in a centralized manner, which requires full communication of system states and control signals on reactive power injections. For the practical implementation considerations of distribution networks, where there is a need for real-time reactive power injections, while the communication infrastructures only support low data rates, the centralized algorithm may not be able to transmit all the state data and control actions through a central operator fast enough~\cite{zhang2014optimal}. In this section, we propose a distributed algorithm to solve \eqref{eqn:main2} by utilizing a trained ICNN for the underlying grid. In the proposed distributed algorithm, we only require communication between neighboring buses and the communication graph is a connected graph~\cite{cady2011robust}. Interestingly, this topology need not to be the same as the physical topology of the distribution network.

\subsection{Problem Formulation}
Again, due to the convexity of ICNN, our problem is essentially following the standard formulation of distributed convex optimization problem~\cite{tsitsiklis1984problems}. This allows us to borrow from a large body of literature on distributed algorithms for convex problems. We are interested in minimizing the sum of their individual convex (potentially non-smooth) objective functions of $N$ interconnected agents:
\begin{equation}
f^{*} :=\min _{q \in \mathpzc{Q}} \quad \sum_{i=1}^{N} f_{i}(\mathbf{q,p}),
\end{equation}
where $\mathpzc{Q}$ is the feasible set defined by the bounds on reactive power injection. Each function $f_i$ is assumed to be pre-trained and known by agent $i$ beforehand. The goal is to solve the minimization problem in a decentralized fashion, where the agents cooperatively
find the optimal reactive power injection without a central coordinator. The coupling can be written out explicitly as
\begin{subequations}
	\label{eqn:distributed2}
	\begin{align}
	\min_{q^{1}, \ldots, q^{N} \in Q} &\sum_{i=1}^{N} f_{i} (\mathbf{q,p}) \\
	\text { s.t. } \quad & q^{1}=q^{2}=\cdots=q^{N}
	\end{align}
\end{subequations}

Problem \eqref{eqn:distributed2} can be solved via dual decomposition~\cite{yu2006dual,chiang2007layering}. Let's denote $q_{j}^{i}$ as node $i$'s estimate of node $j$'s reactive power injection. Then the augmented Lagrangian $L$ is
\begin{equation}
\begin{aligned}
L\left(q^{1}, \ldots, q^{N}, \lambda\right) := \sum_{i=1}^{N} f_{i}\left(q_{1}^{i}, \ldots, q_{N}^{i}\right) +\sum_{i=1}^{N} \sum_{j=1}^{N} \lambda_{i j}\left(q_{i}^{i}-q_{i}^{j}\right)
\end{aligned}
\end{equation}
with corresponding dual function is the infinmum over primal variables:
\begin{equation}
q(\lambda) :=\inf _{q^{1}, \ldots, q^{N} \in Q} L\left(q^{1}, \ldots, q^{N}, \lambda\right)
\end{equation}

Note that $\lambda_{ii}=0$, and we can also rewrite $\sum_{i=1}^{N} \sum_{j=1}^{N} \lambda_{i j} (q_{i}^{i}-q_{i}^{j}) = \sum_{(i, j) \in \mathpzc{E}}\left[
\lambda_{i j}\left(q_{i}^{(i)}-q_{i}^{(j)}\right) +\lambda_{j i}\left(q_{j}^{(i)}-q_{j}^{(j)}\right)\right]$, so there is only $\lambda_{ij}\neq 0$ when there is a communication between neighboring nodes. Then for each node, there is a decomposable subproblem
\begin{equation}
\phi^{i}(\lambda) :=\inf _{q^{i} \in Q} f_{i}\left(q_{1}^{i}, q_{2}^{i}, \ldots, q_{N}^{i}\right)+\sum_{j=1}^{N} \lambda_{i j} q_{i}^{i}-\sum_{j=1}^{N} \lambda_{j i} q_{j}^{i}
\end{equation}

Via gradient ascent, at iteration $t+1$, the update rule for dual variable is
\begin{equation}
\label{eqn:dual_update}
\lambda_{i j}(t+1)=\lambda_{i j}(t)+\alpha_{t}\left(q_{i}^{i}(t)-q_{i}^{j}(t)\right)
\end{equation}
which essentially says that by exchanging Lagrangian multipliers for neighboring nodes, there shall reach  consensus on the dual variables. When all $\lambda_{i j}$ are optimal, $q_i^i$ will be equal to $q_j^i$ for all $(i,j)$ pair, and the duality gap is zero.

So at each iteration, at each agent $i$, by using trained ICNN to parameterize $f_i$, it solves the subproblem simultaneously:
\begin{subequations}
	\label{eqn:distributed}
	\begin{align}
	\min _{q_{i}} \quad & f_{i}\left(q_{1}^{i}, q_{2}^{i}, \ldots, q_{N}^{i}\right)+\sum_{j=1}^{N} \lambda_{i j} q_{i}^{i}-\sum_{j=1}^{N} \lambda_{j i} q_{j}^{i} \\
	\text { s.t. } \quad &\underline {q} \leq q_i \leq \overline{q}
	\end{align}
\end{subequations}
Since the optimization problem is convex, iterating this process~(Algorithm~\ref{algorithm_distributed}) guarantees convergence to the global optimal solution.

\begin{algorithm}
	\label{algorithm_distributed}
	\caption{Distributed Voltage Regulation via ICNN}
	\begin{algorithmic}
		\REQUIRE {Input dataset $\{\mathbf{p,q}\}$}
		\REQUIRE {Number of buses $N$}
		\REQUIRE {Trained ICNN model $h_\theta (\mathbf{p,q})$}
		\WHILE{$\left|q_{i}^{(i)}-q_{i}^{(j)}\right|>\delta$ for any $(i, j) \in \mathpzc{E}$}
		\FOR{$i=1,...,N$}
		\STATE{Solve subproblem  \eqref{eqn:distributed} for $i$ based on $\lambda_{ij}$}
		\STATE{Get value $q_{i}^i$ and $q_{j}^i, \forall j\sim i$}
		\ENDFOR
		\STATE{Update $\lambda_{ij}$ based on $q_{i}^i$ and $q_{i}^j$ using \eqref{eqn:dual_update}}
		\ENDWHILE
		\STATE {Nodal optimal reactive power injection $q_{i}$}
	\end{algorithmic}
\end{algorithm}

\subsection{Smooth Activation Functions}
Even though Algorithm~\ref{algorithm_distributed} is guaranteed to converge, it may converge slowly because the optimization problem is not smooth by the construction described in Section \ref{Sec:ICNN}. Recall that the ReLU activation function is not differentiable at $0$, therefore all of the algorithms we have presented are subgradient algorithms. Even if these subgradients can be computed extremely efficiently through automatic differentiation~\cite{baydin2018automatic}, they can still become a bottleneck for convergence, especially in the distributed setting. Therefore, we can use a smooth version of ReLU as the activation function, which was named as Softplus activation function~\cite{nair2010rectified}. In particular, define $g_t(x)$ as
\begin{equation*}
	g_t(x)=\frac{1}{t} \ln (1+e^{tx}).
\end{equation*}
As $t$ goes to infinity, $g_t$ converges the ReLU function uniformly on $\mathbb{R}$, and it is smooth for any finite $t$. By choosing a large $t$, we can make ICNN smooth from input to output and still maintain a good fit of convex functions.

\section{Simulation Results}
\label{Sec:Experiment}
In this section, we evaluate the proposed voltage regulation scheme on standard distribution networks. Our simulation focus on the IEEE 13-bus and IEEE 123-bus test systems.  Linear models, standard neural networks and the optimal SOCP formulations are used  for comparison.

\begin{figure}[ht]
	\centering
	\includegraphics[width=1.0 \columnwidth]{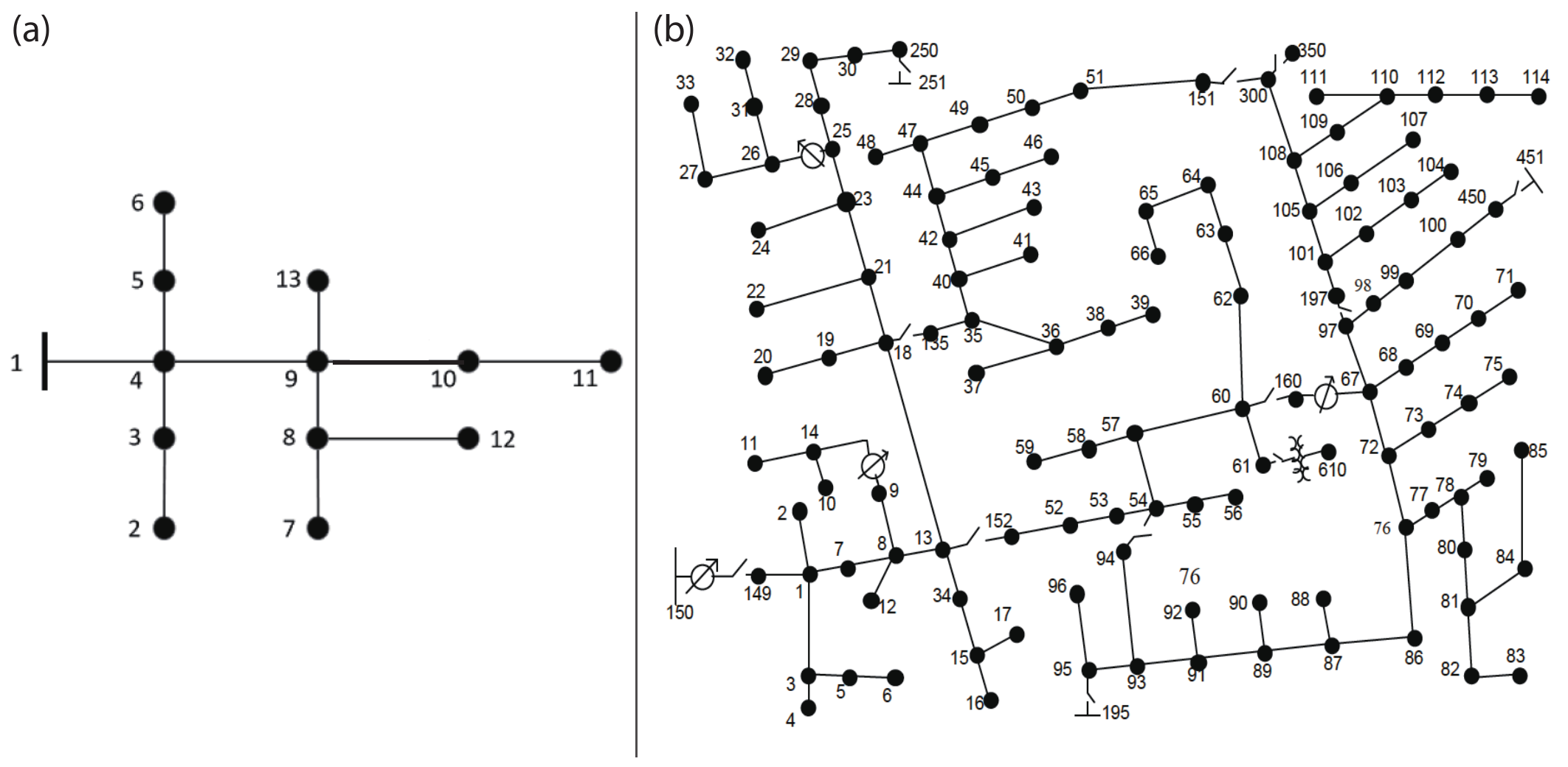}
	\caption{\small Schematic diagram of (a). IEEE 13-bus test feeder and (b). IEEE 123-bus test feeder. Reference buses: 1 and 149.}
	\label{fig:grid}
\end{figure}

\subsection{Simulation Setup}
For both 13-bus and 123-bus system, we use AC power flow model~\eqref{equ:branch_flow} to generate $10,000$ instances of simulation data composed of $\{\mathbf{p},\mathbf{q}, \mathbf{V}\}$.  We assume both the distribution network topology and line parameters are not revealed to the optimization algorithm, except when the optimal SOCP is used as a baseline. We use one year's load data from the University of Washington Seattle campus for training and test dataset generation. We allow plus/minus $20\%$ of reactive power injections at each node as control inputs. We develop three algorithms and compare their performances for two test feeders shown in Fig. \ref{fig:grid}:
\begin{itemize}
	\item \emph{Linear  Model:} We consider using a linear model to fit the unknown dynamics from active and reactive power to the deviations between nodal voltage and the nominal voltage. Such linearized models have been widely used in power systems literature~\cite{zhu2015fast,li2018distribution};
	\item \emph{Neural Networks Model:} We construct standard three-layer and four-layer neural networks for the 13-bus and 123-bus cases, respectively. We tune the parameters of neural networks (e.g., number of neurons, learning rate) and stop the training process once the fitting performance on validation data converges;
	\item \emph{Input Convex Neural Networks:} We keep the number of layers and matrices $W_i, i=1,...,k$ the same dimension as those of neural networks models, but add direct layers $D_i,i=2,...,k$ correspondingly. We constrain network weights $W_{2:k}$ to be non-negative during training.
\end{itemize}

\begin{table*}[h]
	\renewcommand{\arraystretch}{1.4}
	\centering
	\begin{tabular}{>{\centering\arraybackslash}m{15em}|>{\centering\arraybackslash}m{4em}|>{\centering\arraybackslash}m{4em}|>{\centering\arraybackslash}m{4em}|>{\centering\arraybackslash}m{4em}|>{\centering\arraybackslash}m{4.5em}|>{\centering\arraybackslash}m{4em}|>{\centering\arraybackslash}m{4em}|>{\centering\arraybackslash}m{4em}}
		\ChangeRT{1.2pt}
		SImulation Network    & \multicolumn{4}{c|}{IEEE 13-Bus}& \multicolumn{4}{c}{IEEE 123-Bus} \\
		\ChangeRT{0.7pt}
		Model& SOCP & Linear & NN & ICNN & SOCP & Linear&NN&ICNN\\
		\ChangeRT{0.2pt}
		Model Fitting MAE& - &9.93$\%$ &3.45$\%$ &3.86$\%$ &-&12.98$\%$ &3.56$\%$ &4.25$\%$  \\
		\hline
		Regulated voltage out of $3\%$ tolerance &3.46$\%$&8.65$\%$&7.88$\%$&4.71$\%$&3.61$\%$&21.46$\%$&14.04$\%$&7.51$\%$ \\
		Regulated voltage out of $5\%$ tolerance &0.47$\%$&7.89$\%$&6.86$\%$&1.05$\%$&0.72$\%$&19.19$\%$&9.65$\%$&1.64$\%$\\
		\hline
		Computation Time (per instance/s)
		&0.9684&0.2022&0.3137&0.2512&4.041&0.2712&0.6297&0.4302\\
		\ChangeRT{1.2pt}
	\end{tabular}
	\caption{Comparison between SOCP, Linear model, Neural Networks model, and Input Convex Neural Networks model for IEEE 13-bus and IEEE 123-bus systems.}
	\label{table}
\end{table*}

To fit the parameters of neural network models, we use mean squared error as the loss function during training.  To solve the voltage regulation problem \eqref{eqn:voltage_regulation}, we set $\alpha_i,i=1,...,k$ in \eqref{eqn:main2} to be $1$ in our simulation cases. When $\alpha$ is not equal to $1$, we could adapt the optimization problem using weighted sum of voltage deviation correspondingly. Note that we could also flexibly use alternative loss terms or add reactive power costs to the objective function \eqref{subequ:loss}, as long as they are convex functions over reactive power injections. All the implementations are conducted on a MacBook Pro with 2.4GHz Intel Quad Core i5.

\begin{figure}[ht]
	\centering
	\includegraphics[width=0.6 \columnwidth]{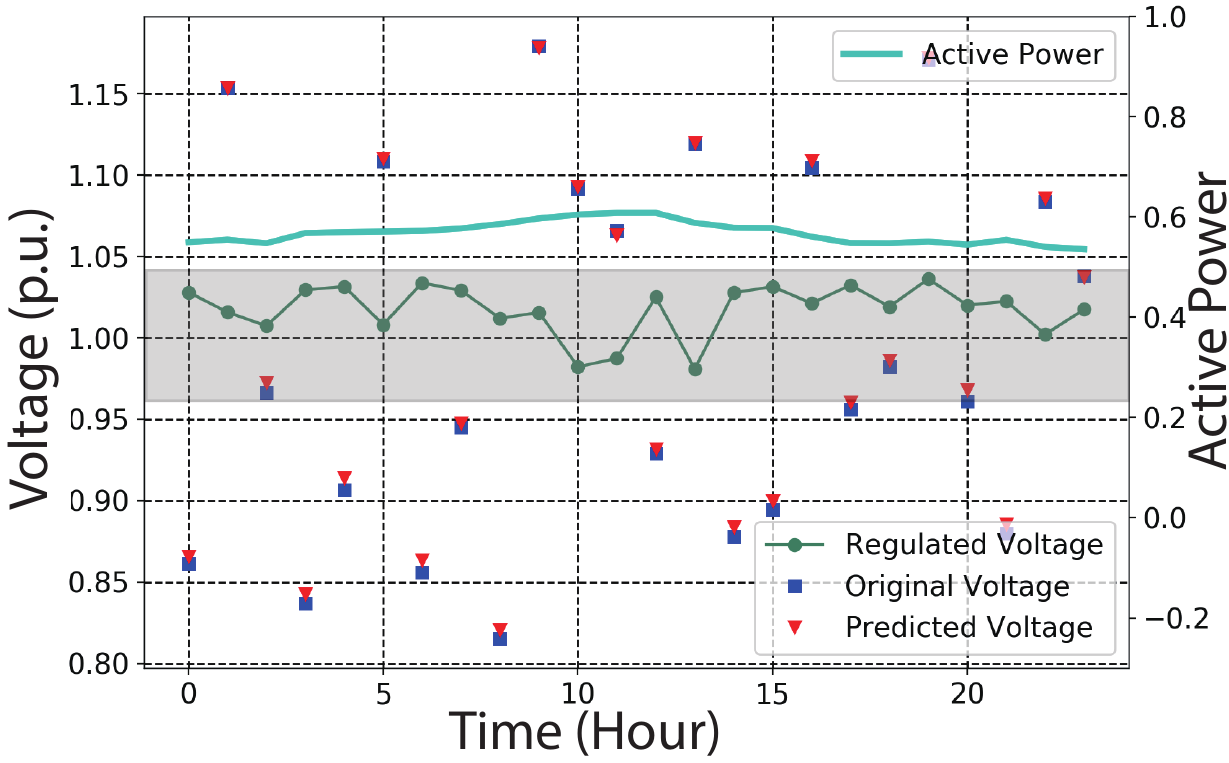}
	\caption{\small Example of voltage regulation over a daily variation for the 13-bus test feeder. The voltage of bus 4 is shown. With ICNN accurately predicting voltages~(red triangle), it could regulate voltage within $4\%$ of nominal values (grey box) under varying load level throughout the day.}
	\label{fig:results}
\end{figure}

To benchmark the performance of the proposed algorithms under unknown topology and parameters, we also follow \cite{farivar2011inverter} to relax $l_{ij}\geq \frac{P_{ij}^2+Q_{ij}^2}{V_i^2}$ in the Dist-flow equations, and use the same validation datasets to solve the resulting convex SOCP. We calculate the optimal reactive power injections along with the resulting voltage profiles. We use CVX to solve the SOCP and linearized models~\cite{grant2014cvx}, and use Tensorflow to set up and optimize over NN and ICNN models~\cite{abadi2016tensorflow}.

\subsection{Estimation Accuracy}

We firstly validate that ICNN can be used as a proxy for power flow equations, and predict the nodal voltage magnitude deviations. By using $8,000$ training instances, the ICNN can predict the voltage deviations on the validation instances accurately. As shown in Table \ref{table}, the mean absolute error~(MAE) of ICNN fitting are smaller than $4.3\%$ in both test systems, which are comparable to $3.45\%$ and $3.56\%$ by using neural networks. This is also illustrated in Fig. \ref{fig:results}, where under different load levels throughout 24 hours, the ICNN can predict all the nodal voltages accurately. More importantly, linear model's fitting performances are over $2$ times worse than the neural networks counterparts. We later show such performances on model tractability and fitting errors would also impact the controller performances.


\begin{figure}[ht]
	\centering
	\includegraphics[width=0.6 \columnwidth]{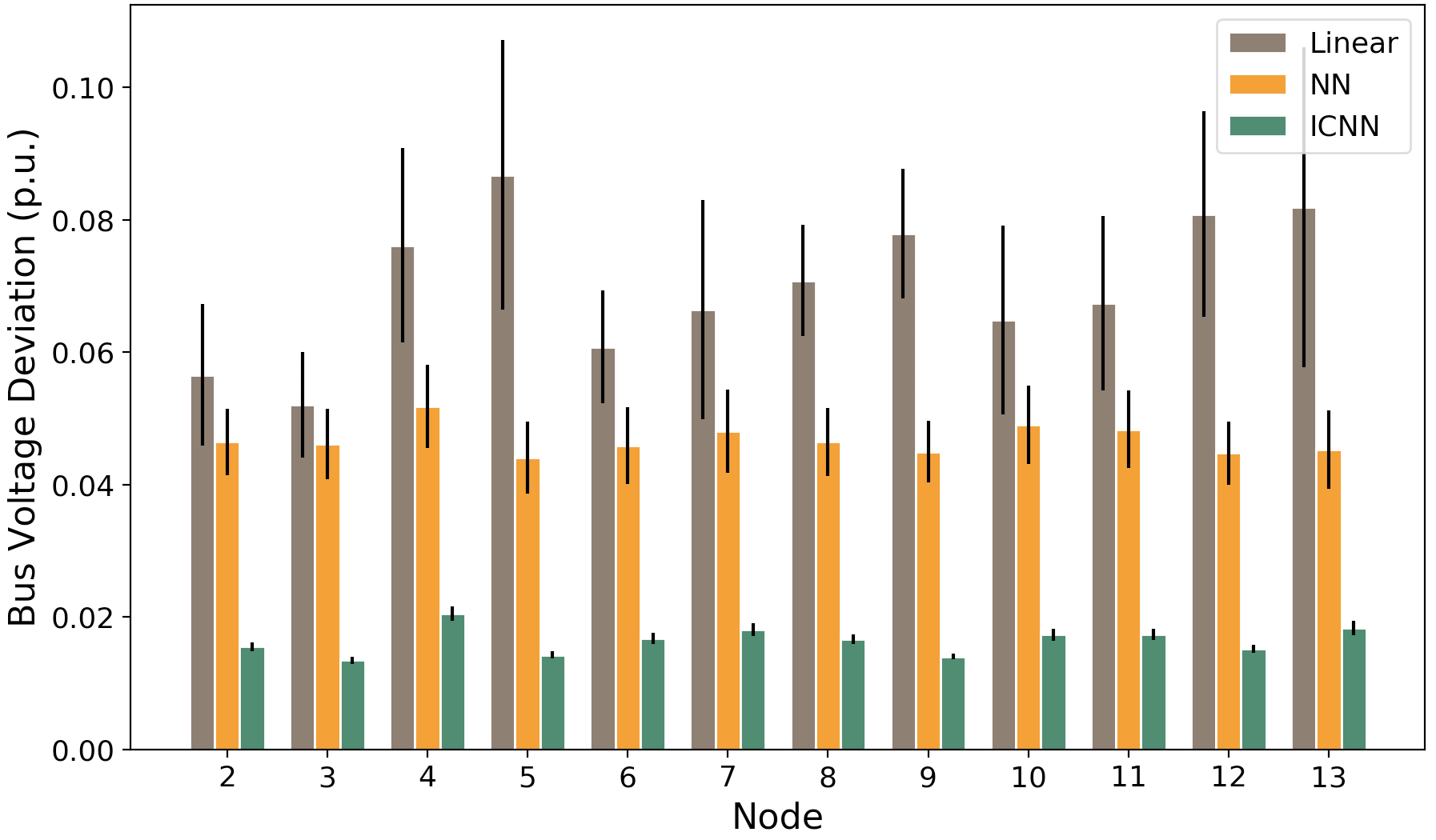}
	\caption{\small Comparisons on nodal voltage deviation bar plots of linear-fitted model, neural network model and input convex neuural network model on IEEE 13-bus system. On average, the mean voltage deviation for ICNN is $4.3$ times smaller than linear model, and $2.7$ times smaller than standard NN model.}
	\label{fig:nodal}
\end{figure}

\begin{figure}[t]
	\centering
	\includegraphics[width=0.6 \columnwidth]{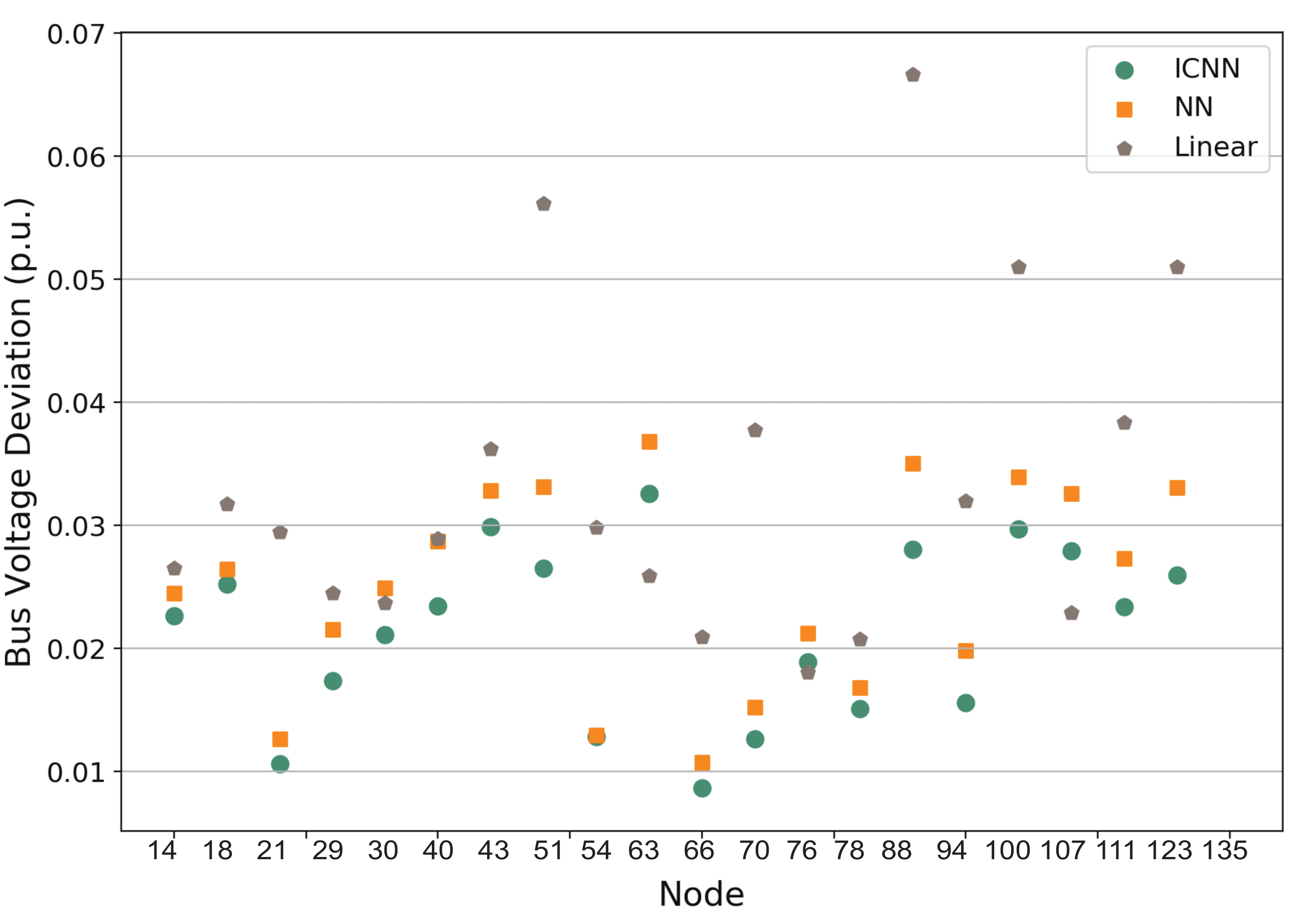}
	\caption{\small Comparisons on 20 randomly selected buses' nodal voltage deviation plots of linear-fitted model, neural network model and input convex neuural network model on IEEE 123-bus system.}
	\label{fig:nodal_voltage}
\end{figure}

\begin{figure*}[t]
	\centering
	\includegraphics[width=1.0 \columnwidth]{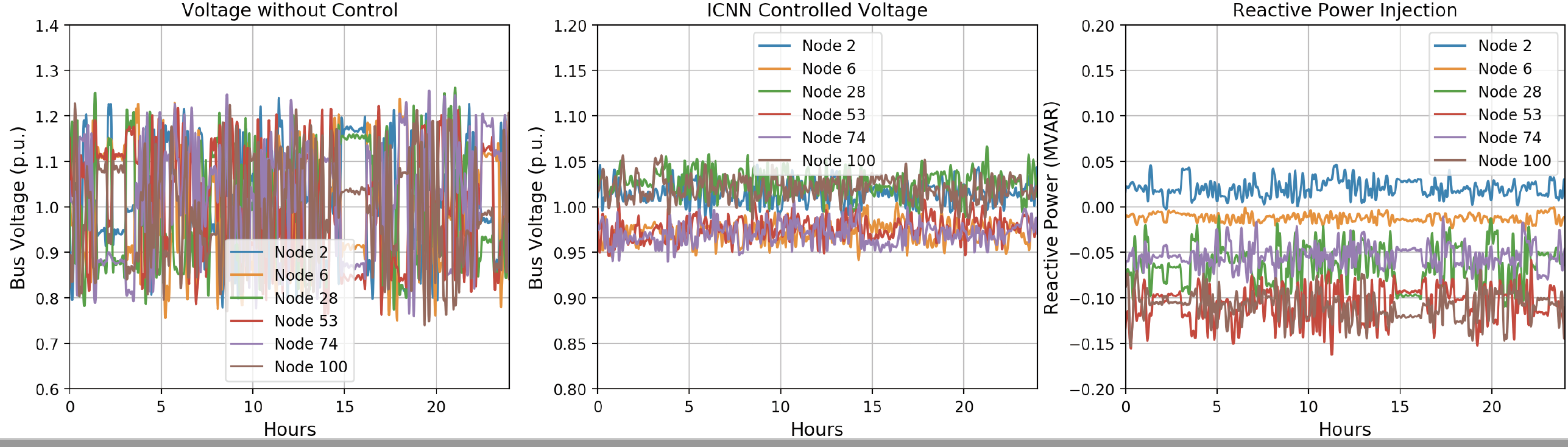}
	\caption{\small Simulation results of proposed distributed algorithm with 24-hour time-varying nodal loads and generations. The left and middle plots show the control performances on nodal voltages, while right plot shows the scale of control actions (reactive power injections).}
	\label{fig:distributed}
\end{figure*}

\subsection{Voltage Regulation Performance}
In Figure \ref{fig:results}, we show the regulated voltage using ICNN in the IEEE 13-bus case. Under this day's load profile, we are able to regulate node 4's voltage magnitude within $\pm 4\%$ per unit with constrained reactive power injections~(Equation \ref{eqn:constraints2}). In Figure \ref{fig:nodal}, we show that the mean and variance on each bus's voltage deviations using three models for the 13-bus feeder. On the one hand, with similar fitting performances, ICNN outperforms the standard neural network in regulating nodal voltages. This is due to the fact that neural networks may have many local minima, and the NN-based controller can not find the optimal reactive power injections. On the other hand, even though linear model provides a easier venue for solving optimization problem, it suffers from inaccurate modeling of the underlying distribution grids, and the regulated bus voltages have greater level of fluctuations. Similar observations also hold in the 123-bus test case, where in Fig. \ref{fig:nodal_voltage} we show the nodal voltage comparison using three models, and voltage regulated by ICNN are constrained to be in a much narrower range. More results on voltage regulation performances are summarized in Table \ref{table}. Under varying load and power generation profiles, ICNN is able to maintain over $98.3\%$ of nodal voltages within $5\%$ deviations from nominal voltages, which are comparable to SOCP solutions. On the contrary, linear fitted models can not scale to larger system, and nearly $20\%$ of voltages are out of $5\%$ tolerance in the 123-bus case.

We also give an analysis on the computation time for each algorithm as shown in last row of Table \ref{table}. Reported results are averaged over 2,000 testing instances. Compared to linear model, optimization based on ICNN generally takes longer time due to the model complexity, but it is still able to find the optimal solutions within the acceptable time range. Note that we are solving the ICNN optimization problem using our own solver, while solving SOCP and linear model using the off-the-shelf CVX solvers. The optimization solution process involving ICNN can be further accelerated with GPU support in the future. More importantly, in the 13-bus case, ICNN-based optimization is faster than SOCP solver, and it scales to 123-bus case with moderate computation time increases compared to SOCP solver. This makes ICNN as a practical modeling and optimization tool for unknown distribution grids. An interesting observation is that it takes longer for NN to find solutions compared to ICNN, partly due to the fact that gradient-based optimizer is stuck in some local minima in normal NN.

\subsection{Distributed Algorithm}
In Figure \ref{fig:distributed} we show the performance of proposed learning and control algorithm under the distributed settings, where the underlying topology is unknown to nodes' controllers, while local communication is allowed for neighboring controllers in the communication graph. We plot the voltage and reactive power injection profiles on 6 randomly selected nodes under 24-hour's varying active power injections. For each node, the voltage magnitude deviation can be controlled effectively using ICNN based controller. As we assume control components at all buses for the distribution grid can supply or consume at most 0.2 MVar reactive power, it is shown in the right plot of Figure \ref{fig:distributed}  that such constraints are satisfied at all buses at all times.

\section{Conclusion}
\label{Sec:Conclusion}
In this paper, we proposed an optimal voltage control framework using neural networks under unknown distribution system topology and parameters. By utilizing the available smart meter data, we design an input convex neural network as an accurate learner for the mapping from power injections to voltage magnitude deviations. We then proposed a distributed algorithm that can be
implemented in distribution grids with a large number of buses using the learned neural network. We demonstrated the performance and efficiency of proposed algorithm in a set of case studies. Future work includes extending proposed  learning and decision-making framework into online settings with dynamic operating scenarios.


\bibliographystyle{IEEEtran}
\bibliography{ref}

\end{document}